\documentclass
[
11pt,
]
{amsart}
\usepackage{
amssymb
}

\newtheorem{theorem}{Theorem}
\newtheorem{proposition}{Proposition}
\newtheorem{lemma}{Lemma}

\theoremstyle{remark}

\theoremstyle{remark}
\newtheorem{remark}{Remark}

\newcommand{\eps}{\varepsilon}
\newcommand{\R}{{\bf R}}
\newcommand{\Id}{\mbox{Id}}

\renewcommand{\r}[1]{(\ref{#1})}
\newcommand{\PDO}{$\Psi$DO}
\newcommand{\be}[1]{\begin{equation}\label{#1}}
\newcommand{\ee}{\end{equation}}

\renewcommand{\d}{\mathrm{d}}

\title[An inverse source problem in optical molecular imaging]{An inverse source problem in optical molecular imaging}

\author[P. Stefanov]{Plamen Stefanov}
\address{Department of Mathematics, Purdue University, West Lafayette, IN 47907}
\thanks{First author partly supported by a NSF FRG Grant No.~0554065}

\author[G. Uhlmann]{Gunther Uhlmann}
\address{Department of Mathematics, University of Washington, Seattle, WA 98195}
\thanks{Second author partly supported by a NSF FRG grant No.~0554571 and a Walker Family Endowed Professorship}

\begin{document}

\begin{abstract} 
We study the direct and an inverse source problem for the radiative transfer equation 
arising in optical molecular imaging.   We show that for  generic  absorption  and scattering coefficients, the direct problem is well-posed and the inverse one is uniquely solvable, with a stability estimate.

\end{abstract}

\maketitle

\section{Introduction}

We consider an inverse source problem arising in optical molecular imaging (OMI) 
which is currently undergoing a rapid expansion. The design
of new biochemical markers that can detect faulty genes and other molecular processes
allows us to detect diseases before macroscopic symptoms appear. This has
been studied extensively in the bioengineering literature. 
See for instance \cite{CBG},
\cite{COSC}, \cite{JRB}.
Unlike higher-energetic markers used in classical nuclear imaging techniques
such as single photon emission computed tomography (SPECT), positron emission
tomography (PET) as well as magnetic resonance imaging (MRI), optical markers emit
relatively low-frequency photons. The objective of OMI is to reconstruct the concentration
of such markers from their radiations measured at the boundary of the domain. 
The radiations in OMI are governed by the equations of radiative transfer and the inverse problem in
OMI is thus an inverse transport source problem, at least once the optical properties of the
underlying medium are known. 
We now describe more precisely the mathematical problem.

We assume that  $\Omega$ is a bounded domain of $\R^n$ with smooth boundary. We will assume also that $\Omega$ is strictly convex. This is not an essential assumption since for the problem that we study, one can always push the boundary away and make it strictly convex, without losing generality.

The radiative transport equation
\be{t1}
\theta\cdot \nabla_x u(x,\theta) +\sigma(x,\theta) u(x,\theta) -\int_{S^{n-1}}k(x,\theta,\theta')  u(x,\theta')\, \d\theta' = f(x),\quad  \quad u|_{\partial_- S\Omega}=0.
\ee 
where the \textit{absorption} $\sigma$ and the \textit{collision kernel} $k$ are functions with a regularity that will be specified below.  The source term $f$ is assumed to depend on $x$ only.

In section~\ref{sec_2} we study the direct problem. We show that for an open and dense set of absorption and scattering coefficients the direct problem \r{t1} is well-posed.  See Theorem~\ref{thm_transport} for details.

The boundary measurements are modeled by
\be{t4}
Xf(x,\theta) = u|_{\partial_+ S\Omega}, \quad (x,\theta)\in \partial_+ S\Omega
\ee
where  $u(x,\theta)$ is a solution of \r{t1}, and $\partial_+ S\Omega$ denotes the points $x\in \partial \Omega$ with
direction $\theta$ pointing outwards.

In section 3 we consider the inverse source problem, that consists
in determining the source term $f$ from measuring $Xf$.  Notice that in the case
$\sigma=k=0$ the linear operator $X$ is the standard X-ray transform and when $k=0$, $X$ is  a weighted
X-ray transform (see section~\ref{sec_2}).

This inverse problem has been considered in several papers in the mathematical
and engineering community \cite{Bal-Tamasan}, \cite{L}, \cite{P}, \cite{Si}, \cite{YSM}.  In particular in \cite{Bal-Tamasan} it is shown
that one can prove uniqueness when  $k=k(x,\theta\cdot\theta')$, and $k$ is small enough in a suitable norm. We show that for the 
absorption and scattering in an dense and open subset we can uniquely determine the source $f$ from the boundary measurements. We also prove a stability estimate. See Theorem~\ref{thm_tr_gen}  for details.

\section{The Direct Problem} \label{sec_2}
Set
\be{t2}
T_0 = \theta\cdot \nabla_x, \quad T_1 = T_0+\sigma, \quad T= T_0+\sigma-K,
\ee
where $\sigma$ is viewed as the operator of multiplication by $\sigma(x,\theta)$, and $K$ is the integral operator in \r{t1}.

Let $u$ solve
\be{t3}
Tu=f, \quad u|_{\partial_- S\Omega}=0.
\ee

As mentioned in the introduction the operator $X$ is the X-ray transform if $\sigma$ and $k$ are both zero
\[
Xf(x,\theta) = If(x,\theta) := \int_{\tau_-(x,\theta)}^0 f(x+t\theta)\, \d t, \quad (x,\theta)\in \partial_+ S\Omega  \quad \mbox{($\sigma=k=0$),}
\]
where $\pm \tau_\pm(x,\theta)\ge0$ are defined by $(x,x+ \tau_\pm(x,\theta)) \in \partial_\pm S\Omega$. 
We will always extend $f$ as $0$ outside $\Omega$, therefore we can assume that we integrate above over $\R$. If $k=0$, then $X$ reduces to the following weighted X-ray transform 
\be{t6}
Xf(x,\theta) = I_\sigma f(x,\theta) := \int E(x+t\theta,\theta) f(x+t\theta)\, \d t, \quad (x,\theta)\in \partial_+S\Omega  \quad \mbox{($k=0$),}
\ee
where 
\be{t7} 
E(x,\theta) = \exp\left(  -\int_0^{\infty}\sigma(x+s\theta,\theta)\, \d s\right).
\ee
If $\sigma>0$ depends on $x$ only, this is known as the attenuated X-ray transform, that is injective, and there is an explicit inversion formula, see \cite{Novikov}, \cite{AB}.

We define the adjoint $X^*$ of $X$ w.r.t.\ the measure $\d\Sigma$ defined above. 
We will view $X$ as a perturbation of $I_\sigma$, and our goal is to show that $X^*X$ is a relatively compact perturbation of $I^*_\sigma I_\sigma$. 

First we will analyze the direct problem. In the next theorem, $f$ is allowed to depend on $\theta$ as well.

\begin{theorem}  \label{thm_transport}
There exists an open and dense set of pairs $(\sigma,k) \in C^2(\bar\Omega\times S^{n-1})\times C^2(\bar\Omega\times S^{n-1}\times S^{n-1})$, including a neighborhood of $(0,0)$, so that for each $(\sigma,k)$ in that set, 

  (a) the direct problem \r{t3}  has a unique solution  $u\in L^2(\Omega\times S^{n-1})$ for any $f\in L^2(\Omega\times S^{n-1})$ depending both on $x$ and $\theta$.

  (b) $X$ extends to a bounded operator  
\[
X : L^2(\Omega\times S^{n-1}) \longrightarrow L^2(\partial_+S\Omega, \, \d\Sigma).
\]
\end{theorem}

\begin{proof}
We start with the analysis of the direct problem \r{t3}. In what follows, let $T_0$, $T_1$ and $T$ denote the operators given by \r{t1} in $L^2(\Omega\times S^{n-1})$ with domain 
\[
D(T_0) = D(T_1)=D(T) = \left\{f\in 
L^2(\Omega\times S^{n-1}); \; \theta\cdot\nabla_x u\in L^2(\Omega\times S^{n-1}), \; u|_{\partial_-S\Omega}=0\right\}.
\]
We will assume here that $f$ depends both on $x$ and $\theta$. Note first that the solution to the problem \r{t3} with $k=0$ is given by $u=T_1^{-1}f$, where 
\be{T1}
[T_1^{-1} f](x,\theta) = \int_{-\infty}^0 \exp\left( -\int_s^0 \sigma(x+\tau\theta,\theta)\, \d \tau  \right)f(x+s\theta,\theta)  \, \d s.
\ee
This follows easily from the fact that $E$ is an integrating factor, i.e., $T_0 = E^{-1}T_1E$. 

Apply $T_1^{-1}$ to both sides of \r{t3} to get
\[
u = T_1^{-1}(Ku+f).  
\]
We therefore  see that  \r{t3} is equivalent to the integral equation
\be{t9n}
(\Id - T_1^{-1}K)u = T_1^{-1}f.
\ee
Therefore, if $\Id - T_1^{-1}K$ is invertible,   \r{t3} is uniquely solvable, and the solution is given by 
\be{t9}
u = T^{-1}f = (\Id - T_1^{-1}K)^{-1}T_1^{-1}f.
\ee
We also get that $T^{-1}$ exists if and only if $(\Id - T_1^{-1}K)^{-1}$ exists. Then
\be{t9xx}
u = T^{-1}f = (\Id - T_1^{-1}K)^{-1}T_1^{-1}f.
\ee
When $f$ depends on $x$ only, set 
\be{t10abc} 
[Jf](x,\theta) := f(x).
\ee
Then
\be{t9x}
u = T^{-1}Jf = (\Id - T_1^{-1}K)^{-1}T_1^{-1}Jf.
\ee

\begin{lemma}  \label{lemma_1}
The operator $KT_1^{-1}J: L^2(\Omega) \to L^2(\Omega\times S^{n-1})$ is compact. 
\end{lemma}

\begin{proof} 
Let first $f$ depend both on $x$ and $\theta$. Then
\be{t12}
\begin{split}
[KT_1^{-1}f](x,\theta) &= \int_{S^{n-1}}k(x,\theta,\theta') \int_{-\infty}^0 \exp\left( -\int_s^0 \sigma(x+\tau\theta',\theta')\, \d \tau  \right)  
f(x+s\theta',\theta')  \, \d s\, \d \theta'\\
 &= \int \frac{\Sigma\left(x,|x-y|,\frac{x-y}{|x-y|}\right)   k \left( x,\theta,\frac{x-y}{|x-y|} \right)}{|x-y|^{n-1}}  f\left(y,  \frac{x-y}{|x-y|} \right)\, \d y,  
\end{split}
\ee
where 
\[
\Sigma(x,s,\theta') = \exp\left( -\int_{-s}^0 \sigma(x+\tau\theta',\theta')\, \d \tau  \right)
\]
(we replaced $s$ by $-s$ and then made the change $x-s\theta'=y$). 

Note that when $f$ depends on $\theta$ as well there are no enough  integrations above to show that $T_1^{-1}K$ is compact. On the other hand, if $f$ depends on $x$ only, i.e., if we have $Jf$ above with such an $f$, we have
\be{KT}  
[KT_1^{-1}J ]f(x,\theta) = \int_\Omega \frac{\Sigma\left(x,|x-y|,\frac{x-y}{|x-y|}\right)   k \left( x,\theta,\frac{x-y}{|x-y|} \right)}{|x-y|^{n-1}}  f(y)\, \d y.
\ee
The integral above is a typical singular operator with a weakly singular kernel, and an additional parameter $\theta$, see \cite{MikhlinP}, \cite{St}. Under the smoothness assumptions on $\sigma$, $k$, it is easy to see that $\partial_\theta KT_1^{-1}$ and $\partial_x
KT_1^{-1}$ are bounded operators, see Proposition~\ref{pr_s} below.  This completes the proof of the lemma.
\end{proof}

As we mentioned above, the arguments above do not prove that $KT_1^{-1}$ is compact on $L^2(\Omega\times S^{n-1})$. On the other hand, its square is compact.

\begin{lemma}  \label{lemma_2}
The operator $T_1^{-1}KT_1^{-1} : L^2(\Omega\times S^{n-1}) \to L^2(\Omega\times S^{n-1})$ is compact. 
\end{lemma}

\begin{proof}
Replace $f\Big(y,\frac{x-y}{|x-y|}\Big)$ in \r{t12} by 
\[
[Kf]\left(y,\frac{x-y}{|x-y|}\right) = \int_{S^{n-1}} k\left(y,\frac{x-y}{|x-y|} ,\theta'\right) f(y,\theta') \, \d \theta'.
\]
Then the compactness follows from the same arguments as in Lemma~\ref{lemma_1}. Indeed, we have
\[
[KT_1^{-1}Kf](x,\theta) = \iint_{\Omega\times S^{n-1}} \frac{\alpha    
\left(x,y,|x-y|,\frac{x-y}{|x-y|},\theta, \theta'\right)  
  }{|x-y|^{n-1}} f(y,\theta') \,\d y \, \d\theta'
\]
with an obvious definition of $\alpha$. In particular, all second  order derivatives of $\alpha$ are bounded. Let $g(x,\theta,\theta')$ be the $y$-integral above, i.e., the r.h.s.\ above becomes $\int g(x,\theta,\theta')\, \d\theta'$.  
Then by Proposition~\ref{pr_s} below,
\[
\int_\Omega |\partial_x g(x,\theta,\theta')|^2\, \d x \le C\int_\Omega |f(x,\theta')|^2\, \d x
\]
for any $\theta$, $\theta'$. In particular,
\[
\iint_{\Omega\times S^{n-1}} |\partial_x g(x,\theta,\theta')|^2\, \d x \, \d \theta' \le C\|f\|^2_{L^2}. 
\]
Then
\[\begin{split}
\left\|\partial_x KT_1^{-1}Kf\right\|^2  &= \iint_{\Omega\times S^{n-1}}    \Big|\int_{S^{n-1}}\partial_x g(x,\theta,\theta')\, \d\theta'\Big|^2\d x\, \d\theta\\
&\le C\iint_{\Omega\times S^{n-1}}    \int_{S^{n-1}}|\partial_x g(x,\theta,\theta')|^2\, \d\theta' \d x\, \d\theta\\
&\le C'\|f\|^2_{L^2}.
\end{split}
\]
It is easy to see that $\partial_\theta KT_1^{-1}Kf\in L^2$ as well. This, and the estimate above, imply the compactness of $KT_1^{-1}K$. 
\end{proof}

We proceed with the proof of part (a) of the theorem. We are looking for $k$ so that $T^{-1}$ exists. 
Consider
\[
A(\lambda) = \left( \Id -\big(\lambda KT_1^{-1}\big)^2\right)^{-1}
\]
in $L^2(\Omega\times S^{n-1})$. 
The operator $\big(KT_1^{-1}\big)^2$ is compact, and for $\lambda=0$, the resolvent above exists. 
By the analytic Fredholm theorem \cite{Reed-Simon1}, $A(\lambda)$ is a meromorphic family of bounded operators. In particular, it exists for all but a discrete set of $\lambda$'s. Thus for the  those $\lambda$'s, the resolvent $\left( \Id -\lambda KT_1^{-1}\right)^{-1}$ exists and is given by
\be{tl}
\left( \Id -\lambda KT_1^{-1}\right)^{-1} = \left( \Id +\lambda KT_1^{-1}\right) A(\lambda). 
\ee
Indeed, it is obvious that the operator on the r.h.s.\ above is a right inverse to $\Id -\lambda KT_1^{-1}$. For $|\lambda|\ll1$, one can use Neumann series to show that it is left inverse as well. One can use analytic continuation around the poles to show that this remains true   for all $\lambda$ that are not poles.

By \r{t9}, then $T^{-1}$ exists for such $\lambda$'s and $k$ replaced by $\lambda k$. In particular, this shows that the set of such $(k,\sigma)$ is dense. 
Standard perturbation arguments show that the set of $k$'s for which $\Id -\lambda KT_1^{-1}$ is invertible, is open in $C^0$ for a fixed $\sigma$; and the set of pairs $(\sigma,k)\in C^0\times C^0$ with the same property is open, too. Since we just showed that it is dense  as well in $C^0\times C^0$, this completes the proof of (a). 

\bigskip

We proceed with the proof of (b). For  $X$ we get, see \r{t9},
\be{t10x}
Xf = R_+T^{-1}f   = R_+(\Id - T_1^{-1}K)^{-1}T_1^{-1}f,
\ee
where
\[
R_+h = h|_{\partial_+S\Omega}. 
\]
If $f$ depends on $x$ only, then
\be{t10}
Xf = R_+T^{-1}Jf   = R_+(\Id - T_1^{-1}K)^{-1}T_1^{-1}Jf.
\ee 
Notice first that
\be{t10a}
(\Id - T_1^{-1}K)^{-1}T_1^{-1} = T_1^{-1}(\Id - KT_1^{-1})^{-1},
\ee
and in particular, the resolvent on the left exists if and only if the resolvent in the r.h.s.\ does. We therefore have
\[
Xf = R_+ T_1^{-1}(\Id - KT_1^{-1})^{-1}f.
\]
To prove (b), it is enough to show that
\[
R_+ T_1^{-1} : L^2(\Omega\times S^{n-1}) \longrightarrow L^2(\partial_+ S\Omega,\,\d\Sigma)
\]
is bounded. A straightforward computation (see also \cite{Choulli-Stefanov}) shows that
\[
\int_{\partial_+S\Omega}\int_{\tau_-(x,\theta)}^0 f(x-t\theta,\theta)\, \d t\,\d\Sigma  = \int_{\Omega\times S^{n-1}}f(x,\theta)\, \d x\, \d\theta
\]
for any $f\in L^1(\Omega\times S^{n-1})$. Therefore,
\[
\begin{split}
\|R_+T_1^{-1} f\|^2_{  L^2(\partial_+ S\Omega,\,\d\Sigma)} &= 
\int_{\partial_+ S\Omega} |R_+T_1^{-1}f(x,\theta)|^2\d\Sigma 
\le  \int_{\partial_+ S\Omega}\left| \int_{\tau_-(x,\theta)}^0
f(x+t\theta,\theta)\, \d t\right|^2\d\Sigma\\
&\le \int_{\partial_+ S\Omega}\left( |\tau_-(x,\theta)| \int_{\tau_-(x,\theta)}^0
|f(x+t\theta,\theta)|^2\, \d t\right)\d\Sigma\\ 
&\le \mbox{diam}(\Omega) \, \|f\|^2_{L^2(\Omega\times S^{n-1})}.
\end{split}
\]
This completes the proof of Theorem~\ref{thm_transport}. 
\end{proof}

\section{The Inverse Source Problem} \label{sec_3}

In this section we consider the inverse source problem. The next theorem shows that for generic $(\sigma,k)$ there is uniqueness and stability. 
As mentioned in the introduction a similar result has been proven in \cite{Bal-Tamasan} in the case where $k=k(x,\theta\cdot\theta')$, and $k$ is small enough in a suitable norm. 

Fix a bounded domain $\Omega_1$ so that $\Omega_1\Supset\Omega$. Extend $(\sigma,k)$ with regularity as below to functions in $\Omega_1$ with the same regularity.  We chose and fix that extension as  a continuous operator in those spaces. Define the operator  $X_1: L^2(\Omega_1) \to L^2(\partial_+SM_1)$ in the same way as $X$. 
We will be interested in the restriction of $X_1$ to functions  $f$ supported in $\bar \Omega$. We always extend such $f$ as zero to $\Omega_1$. This corresponds to taking measurements on $\partial\Omega_1$ instead of $\partial\Omega$. 

\begin{theorem}  \label{thm_tr_gen} 
There exists an open and dense set of pairs 
\be{t9c}
(\sigma,k) \in C^2(\bar \Omega\times S^{n-1})\times  C^2\left(\bar \Omega_x\times S^{n-1}_{\theta'}; \; C^{n+1}(S^{n-1}_\theta)\right),
\ee
including a neighborhood of $(0,0)$, so that for each $(\sigma,k)$ in that set, the conclusions of Theorem~\ref{thm_transport} hold in $\Omega_1$, and 

  (a) the map $X_1$ is injective on $L^2(\Omega)$,
  
   (b) the following stability estimate holds
\be{t8}
\|f\|_{L^2(\Omega)}  \le C\|X^*Xf\|_{H^1(\Omega_1)}, \quad \forall f\in L^2(\Omega),
\ee
with a constant $C>0$ locally uniform in $(\sigma,k)$.
\end{theorem}

\begin{remark}
The smoothness requirement on $k$ can be reduced to $k\in C^2$ if $k$ is of a special form, like $k=\Theta(\theta)\kappa(x,\theta')$ or a finite sum of such, see \r{t20}, \r{t20a}. 
\end{remark}

>From now on, we will drop the subscript $1$, and all operators below are as defined before but in the domain $\Omega_1$. We assume that $(\sigma,k)$ are such that $T^{-1}$ exists. 
We assume now that $X$ is applied to $f$ that depends on $x$ only. For now, it is not important that $f$ is supported in $\bar\Omega$; that will be needed in \r{21q} and after that; so we apply $X$ to functions in $L^2(\Omega_1)$. 
By \r{t10},
\be{tL}
X = I_\sigma + L,  \quad L := R_+\left(-\Id+(\Id-T_1^{-1}K)^{-1}\right) T_1^{-1}J,
\ee
see also \r{t6}. 
Then
\be{t13}
X^*X = I_\sigma^* I_\sigma +\mathcal{L}, \quad  \mathcal{L} :=  I^*_\sigma L + L^* I_\sigma  + L^*L. 
\ee

In our analysis, we will apply a parametrix of $I_\sigma^* I_\sigma$ to $X^*X$. That parametrix is a first order operator. For this reason, we study $\partial_x I^*_\sigma L$.

\begin{lemma}   \label{lemma_c} The operators
\[
\partial_x I^*_\sigma L, \quad          \partial_x L^* I_\sigma,     \quad \partial_x L^* L 
\]
are compact as operators mapping $L^2(\Omega_1)$ into $L^2(\Omega_1)$.
\end{lemma}

\begin{proof}
To analyze $I^*_\sigma L $, note that $L$ also admits the following representation
\be{t14}
L = R_+T_1^{-1}K T_1^{-1}\left( \Id-KT_1^{-1} \right)^{-1}J.
\ee
We need to study $I^*_\sigma  R_+ T_1^{-1}KT_1^{-1}h$, where $h=h(x,\theta)$. 
Notice first that
\be{t15}
[I^*_\sigma h](x) = \int_{S^{n-1}} \bar E(x,\theta) h^\sharp (x,\theta) \, \d \theta,
\ee
where $\bar E$ denotes complex conjugate, and $h^\sharp$ is the extension of $h\in C(\partial_+S\Omega_1)$ as a constant along the lines originating from $x$ in the direction $-\theta$, see e.g., \cite[sec.~4]{FSU}. In other words, $h^\sharp(x,\theta) = h(x+\tau_+(x,\theta),\theta)$. Next, $R_+T_1^{-1}h$ looks just like $I_\sigma$, see \r{t6} but with $f$ there depending on $\theta$ as well. Therefore,
\[
[I^*_\sigma  R_+ T_1^{-1}]g(x) = \int_{S^{n-1}}  \bar E(x,\theta) \left[\int_{-\infty}^0 E(x+t\theta,\theta)g(x+ t\theta,\theta) \, \d t\right]^\sharp  \, \d \theta.
\]
This yields (see \cite{FSU} again):
\be{t19c}
\begin{split}
[I^*_\sigma  R_+ T_1^{-1}g](x) &= \iint_{S^{n-1}}  \bar E(x,\theta) (E g)(x+t\theta,\theta) \, \d \theta\, \d t\quad\\
&= 2\int_{\Omega_1}\frac{ \Big[\bar E \left(x, \frac{y-x}{|y-x|} \right )(Eg)\left(y, \frac{y-x}{|y-x|} \right)\Big]_{\rm even}}{|y-x|^{n-1}}\, \d y,
\end{split}
\ee
where $F_{\rm even}(x,\theta)$ is the even part of $F$ as a function of $\theta$. If we set $g=KT_1^{-1}h$,  that would give us $I^*_\sigma  R_+ T_1^{-1}KT_1^{-1}h$. 

Instead of assuming \r{t9c}, we will make the following weaker assumption at this point: $k$ can be written as the infinite sum 
\be{t20}
k(x,\theta,\theta')  = \sum_{j=1}^\infty \Theta_j(\theta)\kappa_j(x,\theta')
\ee
with some functions $\Theta_j$ and $\kappa_j$  so that 
\be{t20a}
\sum_{j=1}^\infty  \|\Theta_j\|_{H^1(S^{n-1})} \|\kappa_j\|_{L^\infty(\Omega_1\times S^{n-1})} <\infty.
\ee
One such way to do this is to choose $\Theta_j$ to be the spherical harmonics $Y_j$; then $\kappa_j$ are the corresponding  Fourier coefficients. Then $\|Y_j\|_{H^1(S^{n-1})}\le C(1+\lambda_j)$, where $\lambda_j^2$ are the eigenvalues of the positive Laplacian on $S^{n-1}$. Since $\lambda_j= O(j^{1/(n-1)})$, for the uniform convergence of \r{t20} it is enough to have $\|\kappa_j\|_{L^\infty}\le C (1+\lambda_j)^{-n-\eps}$, $\eps>0$. This would be guaranteed if $k\in L^\infty\left(\Omega_1\times S^{n-1}_{\theta'}; \; C^{n+1}_\theta(S^{n-1})\right)$ by standard integration by parts arguments. Therefore, the hypothesis \r{t9c} of the theorem implies \r{t20}, \r{t20a}.

Under this assumption, for $K_jT_1^{-1}h$, where $K_j$ has kennel $\Theta_j\kappa_j$, we have, see \r{t12},
\be{t21}
\begin{split}
[K_jT_1^{-1}h](x,\theta) 
 &= \Theta_j(\theta) [B_j h](x),\\ 
B_j h(x) :&= \int_{\Omega_1} \frac{\Sigma\left(x,|x-y|,\frac{x-y}{|x-y|}\right)   \kappa_j \left( x,\frac{x-y}{|x-y|} \right)}{|x-y|^{n-1}}  h\left(y,  \frac{x-y}{|x-y|} \right)\, \d y. 
\end{split}
\ee
We claim now that 
$B_j\left(\Id-K T_1^{-1}\right)^{-1}J : L^2(\Omega_1) \to  L^2(\Omega_1)$ is compact. 
We have
\[
\left(\Id-K T_1^{-1}\right)^{-1}J  = J + \left(\Id-K T_1^{-1}\right)^{-1} K T_1^{-1}J.
\]
By Lemma~\ref{lemma_1}, the second term on the right is compact. Therefore, it remains to show that $B_jJ$ is compact. Observe that $B_jJh$  is given by \r{t21} with $h=h(x)$. The compactness then follows from  Proposition~\ref{pr_s}, assuming that $\kappa_j\in C^2$. On the other hand, $B_jJ$ is compact under the assumption that $\kappa_j\in L^\infty$ only, by \cite[Theorem~VII.3.3]{MikhlinP}. Moreover, its norm is bounded by $C\|\kappa_j\|_{L^\infty}$. 

We can write now
\be{t22}
\begin{split}
\partial_x I_\sigma^*L  &= \partial_x 
I_\sigma^*R_+T_1^{-1}KT_1^{-1}\left(\Id-K T_1^{-1}\right)^{-1}J\\
 &= \sum_{j=1}^\infty \left[ \partial_x I_\sigma^*R_+T_1^{-1} \Theta_j J \right]   \left[B_j\left(\Id-K T_1^{-1}\right)^{-1}J\right].
\end{split}
\ee
We notice first that $\partial_x I_\sigma^*R_+T_1^{-1} \Theta_j J : L^2(\Omega_1) \to L^2(\Omega_1)$ is bounded by Proposition~\ref{pr_s}(b), compare to \r{t19c}, with a norm bounded by $C\|\sigma\|_{C^2}\|\Theta_j\|_{H^1}$.  The operator in the second brackets is compact by the claim above. Therefore, each summand in the r.h.s.\ of \r{t22} is a compact operator with a norm 
 not exceeding $C\|\Theta_j\|_{H^1} \|\kappa_j\|_{L^\infty}$, where $C$ depends on $\sigma$ as well.  
Then the series in \r{t22} converges uniformly by \r{t20a}. 
Under this condition,  $\partial_x I_\sigma^*L$ is compact.

To analyze $\partial_x L^* L$, we will follow the proof above. It is enough to show that $\partial_x L^*R_+T_1^{-1} \Theta_j J : L^2(\Omega) \to L^2(\Omega_1)$ is bounded.  
We have, see \r{t14},
\be{t22f}
\begin{split}
\partial_x L^*R_+T_1^{-1} \Theta_j J   &= 
\partial_x \Big( 
R_+T_1^{-1}\left( \Id-KT_1^{-1} \right)^{-1}K T_1^{-1}J
\Big)^* R_+T_1^{-1} \Theta_j J \\
   & = \partial_x \left( K T_1^{-1}J\right)^* \Big( 
R_+T_1^{-1}\left( \Id-KT_1^{-1} \right)^{-1}
\Big)^* R_+T_1^{-1} \Theta_j J.
\end{split}
\ee
Since $R_+T_1^{-1}$ is bounded, it remains to show that the operator $\partial_x \left( K T_1^{-1}J\right)^* : L^2(\Omega_1\times S^{n-1}) \to L^2(\Omega)$ is bounded, as well. The  kernel of the latter is, see \r{KT},
\[
(x,(y,\theta)) \quad \longmapsto \quad \partial_x \frac{\Sigma\left(y,|y-x|,\frac{y-x}{|y-x|}\right)   k \left( y,\theta,\frac{y-x}{|y-x|} \right)}{|y-x|^{n-1}}.
\]
Then the boundedness of  $\partial_x \left( K T_1^{-1}J\right)^*$ then follows as in Lemma~\ref{lemma_2}.

Finally, the fact that $\partial_xL^*I_\sigma$ is bounded follows from the proof for $\partial_xL^*L$. Indeed, $\partial_xL^*I_\sigma = \partial_x L^*R_+T_1^{-1} J $, compare with \r{t22f}, where we can set $\Theta_j=1$. 

This completes the proof of Lemma~\ref{lemma_c}. 
\end{proof}

\begin{proof}[Proof of Theorem~\ref{thm_tr_gen}]

We return to the analysis of the operator $X^*X$, see \r{t13}. We showed in Lemma~\ref{lemma_c} that, up to a relative compact operator, $X^*X$ coincides with $I_\sigma^* I_\sigma$. Assume that $\sigma$ and $k$ are $C^\infty$. Let $Q$ be a parametrix (of order $1$) to the elliptic \PDO\ $I_\sigma^* I_\sigma$ in $\Omega_1$. We restrict the image of $Q$ to $L^2(\Omega)$, i.e., we view $Q$ as an operator $Q : H^1(\Omega_1)\to L^2(\Omega)$. 
Then for any $f$ supported in $\bar\Omega$, we have
\be{21q}
QI_\sigma^* I_\sigma f = f + K_1 f,
\ee
where $K_1$ is of order $-1$ near $\Omega$. Apply $Q$ to $X^*X$ to get
\be{t22a}
QX^*Xf = f +K_2f, \quad K_2:=  K_1 + Q\mathcal{L}.
\ee
Then $K_2: L^2(\Omega) \to L^2(\Omega)$ is compact. We get that the problem of inverting $X^*X$ is reduced to a Fredholm equation. We will show, that it is generically solvable, as in the theorem.

We show first  that the set of pairs  for which $X$ is injective is dense. 

By the results in \cite{FSU}, see  Theorems~1 and 2 there,  if $\sigma$ is real analytic in a $\bar\Omega_1$, then   $I_\sigma$ is injective, and therefore $I_\sigma^* I_\sigma$, is injective as well. Moreover, a small  $C^2(\bar\Omega)$, perturbation preserves that property.   Actually, the remark after the \cite[Theorem~2]{FSU} shows that this is true even for small enough $C^1$ perturbations. Fix $\sigma$ real analytic in $\bar\Omega_1$.   Fix $k$ as well so that $(\sigma,k)$ belongs to the generic set in Theorem~\ref{thm_transport}, related to $\Omega_1$, and the regularity assumption \r{t9c} is satisfied. That can be done for an open dense set of $k$'s by the proof of Theorem~\ref{thm_transport}. 
Consider $X$ related to $(\sigma,\lambda k)$ with $\lambda$ belonging to some complex neighborhood $\mathcal{C}$ of  $[0,1]$. The operator $K_2$ in \r{t22a}  depends meromorphically on $\lambda\in  \mathcal{C}$. Indeed, $K_1$ is related to $(\sigma,k)$ (i.e., to $\lambda=1$), and is therefore independent of $\lambda$. The parametrix $Q$ is also independent of $\lambda$. The analysis above shows that $\mathcal{L}$ is a meromorphic function of $\lambda$ because $L$ has that property, see \r{tl} and \r{tL}. For $\lambda=0$, we have $\mathcal{L}=0$, and then $K_2=K_1$. By adding a finite rank operator to $Q$, we can arrange that $\Id+K_1$, see \r{21q}, is injective, see also the proof of Proposition~4 in \cite{SU-JAMS}. We can then apply the analytic Fredholm theorem again in  $\mathcal{C}$ with the poles of $(\Id-\lambda K)^{-1}T_1^{-1}$ removed. The latter is a connected set, containing $\lambda=0$ and $\lambda=1$. The analytic Fredholm theorem then implies that $QX^*X$ is invertible for all $\lambda$ in that set with the possible exception of a discrete set. In particular, there are $\lambda$'s as close to $\lambda=1$ as needed with that property. For those $\lambda$'s, $X^*X$ and $X$ are injective as well. This shows that there is a dense set of pairs $(\sigma,k)$ in the space \r{t9c} so that $X$ is injective. Lets us call that set $\mathcal{U}$. 

We show next that for $(\sigma,k)$ in some neighborhood of $\mathcal{U}$, $X$ is still injective. 

Let $(k,\sigma)\in \mathcal{U}$. Then   $X: L^2(\Omega)\to L^2(\partial\Omega_1,\d\Sigma)$ is  injective . Then $X^*X : L^2(\Omega) \to H^1(\Omega_1)$ is injective as well, as an integration by parts shows. By adding a finite rank operator to $Q$, we can arrange that $\Id+K_1$, see \r{21q}, is injective, as above.  Then $\Id+K_1$ is invertible on $L^2(\Omega)$, and we deduce that \r{t8} holds.

The analysis above  implies that the norm  $\|X^*X\|_{L^2(\Omega)\to H^1(\Omega_1)}$ depends continuously on $(\sigma,k)$ as in \r{t9c}. Therefore, we can perturb $(\sigma,k)$, and \r{t8} would remain true because the perturbation of the r.h.s.\ will be absorbed by the l.h.s. On the other hand, injectivity of $X^*X$ implies injectivity of $X$.

This proves that the set of pairs $(\sigma,k)$, for which $X$ is injective, is open subset of the (generic set) of pairs, for which the direct problem is guaranteed to be uniquely solvable by Theorem~\ref{thm_transport}. Moreover, \r{t8} with $C$ locally uniform.

This completes the proof of Theorem~\ref{thm_tr_gen}. 
\end{proof}

In the proof of the theorem, we used the following proposition about singular operators. 

\begin{proposition}  \label{pr_s}
Let $A$ be the operator
\[
[Af](x) = \int \frac{ \alpha\Big(x,y,|x-y|,\frac{x-y}{|x-y|} \Big) }{|x-y|^{n-1}} f(y)\, \d y
\]
with $\alpha(x,y,r,\theta)$ compactly supported in $x, y$. Then

(a) If $\alpha\in C^2$, then $A : L^2 \to H^1$ is continuous with a norm not exceeding $C\|\alpha\|_{C^2}$. 

(b) Let $\alpha(x,y,r,\theta) = \alpha'(x,y,r,\theta)\phi(\theta)$. Then $\|A\|_{L^2\to H^1}\le C\|\alpha'\|_{C^2} \|\phi\|_{H^1(S^{n-1})}$.
\end{proposition}

\begin{proof}
We recall some facts about the Calder\'on-Zygmund theory of singular operators, see \cite{MikhlinP}. First,
if $K$ is an integral operator with singular kernel $k(x,y) = \phi(x,\theta)r^{-n}$, where $\theta = (x-y)/|x-y|$, $r=|x-y|$, and if the ``characteristic''  $\phi$ has a mean value $0$ as a function of $\theta$, for any $x$, then $K$ is a well defined operator on test functions, where the integral has to be understood in the principle value sense. Moreover, $K$ extends to a bounded operator to $L^2$ with a norm not exceeding $C\sup_x \|\phi(x,\cdot)\|_{L^2(S^{n-1})}$, see \cite[Theorem~XI.3.1]{MikhlinP}. The characteristic $\phi$ does not need to have zero mean value in $\theta$ but then the integral has to be considered as a convolution in distribution sense. The latter is well defined because the Fourier transform of the kernel w.r.t.\ the variable $z=r\theta$ is homogeneous of order $0$, thus bounded.

Also, see \cite[Theorem~XI.11.1]{MikhlinP},  if $B$ is an operator with a weakly singular kernel $\psi(x,\theta)r^{-n+1}$, then $\partial_xB$ is an integral operator with singular kernel $\partial_x[\beta(x,\theta)r^{-n+1}]$. The latter, up to a weakly singular operator,  has a singular kernel of the type $\phi r^{-n}$, and the integration  is  again understood in the principle value sense, see the next paragraph. In particular, the zero mean value condition is automatically satisfied. 

In our case, $\beta=\alpha$ depends on $y$ and $r$ as well. Assume first that it does not, i.e., $B$ is as above. Extend $\beta$ as a homogeneous function of $\theta$ of order $0$ near $S^{n-1}$. Then
\begin{align}  \nonumber
\partial_{x_i} \frac{\beta(x,\theta)}{r^{n-1}} &= (1-n) \frac{\theta_i}{r^n}\beta + \sum_j\frac{\partial\beta/\partial\theta_j}{r^{n-1}} \frac{\partial\theta_j}{\partial x_i}
+ \frac{ \beta_{x_i}(x,\theta)}{r^{n-1}}\\  \nonumber
 & = (1-n) \frac{\theta_i}{r^n}\beta + \sum_j\frac{\partial\beta/\partial\theta_j}{r^n}  \left(   \delta_{ij}-\theta_i\theta_j   \right)
+ \frac{ \beta_{x_i}(x,\theta)}{r^{n-1}}\\
   &= \frac{(1-n) \theta_i\beta +\partial\beta/\partial\theta_i}{r^n}   
+ \frac{ \beta_{x_i}(x,\theta)}{r^{n-1}}   \label{t23}
\end{align}
We used the fact that $\sum_j\theta_j \partial\beta/\partial\theta_j=0$ because $\beta$ is homogeneous of order $0$ in $\theta$. It is not hard to show that the ``characteristic'' $\phi(x,\theta)= (1-n) \theta_i\beta +\partial\beta/\partial\theta_i$ has zero mean over $S^{n-1}_\theta$, see \cite[p.~243]{MikhlinP}. In this particular case ($\alpha(x,y,\theta)=\beta(x,\theta)$, independent of $y$, $r$), statement (a) can be proven as follows. Choose a finite atlas of charts for $S^{n-1}$ so that for each chart, $n-1$ of the $\theta$ coordinates (that we keep fixed in $\R^n$) can be chosen as local coordinates. By rearranging the $x$, and respectively, the $\theta$ coordinates, in each fixed chart, we can assume that they are $\theta'=(\theta_1,\dots, \theta_{n-1})$. Then $\partial\beta/\partial\theta_n = -\sum_{i=1}^{n-1}\partial\beta/\partial \theta_i$. Then in \r{t23}, we have derivatives of $\beta$ w.r.t.\ $\theta'$ (and $x$)  with smooth coefficients. The contribution of the first term then can be estimated by the Calder\'on-Zygmund theorem. The second term is a kernel of a weakly singular operator. The following criterion can be applied to it: If $K$ has an integral kernel $k(x,y)$ with the property
\be{t23k}
\sup_x \int|k(x,y)|\d x \le M, \quad \sup_y \int|k(x,y)|\d y \le M,
\ee
then $K$ is bounded in $L^2$ with a norm not exceeding $M$ \cite[Prop.~A.5.1]{Taylor-book1}. 

This proves (a) for $\alpha=\beta$. 

To replace $\beta(x,\theta)$ above by $\alpha(x,y,\theta)$, write
\[
\alpha(x,y,r,\theta) = \alpha(x,x,0,\theta) +r\gamma(x,y,r,\theta). 
\]

To prove (b), write first as above,
\be{t24}
\alpha(x,y,r,\theta) = \beta'(x,\theta)\phi(\theta) +r\gamma(x,y,r,\theta)\phi(\theta), \quad \beta'(x,\theta) := \alpha_1(x,x,0,\theta),
\ee
where $\gamma\in C^1$. Notice then that in \r{t23}, with $\beta=\beta'\phi$, we have
\[
(1-n) \theta_i\beta +\partial\beta/\partial\theta_i = (1-n) \theta_i\beta'\phi +\phi\partial\beta'/\partial\theta_i +\beta'\partial\phi/ \partial\theta_i.
\]
Choosing local coordinates as above, and applying the Calder\'on-Zygmund theorem again, we get that the first term above contributes a singular operator with a norm not exceeding $\|\alpha_1\|_{C^1}\|\phi\|_{H^1}$. The second term $r\gamma$ generates an operator with a kernel  $\gamma(x,y,r,\theta)\phi(\theta)r^{-n+2}$. Differentiate w.r.t.\ $x$, and we still get a weakly singular operator which norm can be estimated as in \r{t23k} to give a norm not exceeding $\|\gamma_1\|_{C^1}\|\phi\|_{H^1}$. 
\end{proof}

\begin{remark}
The only second order derivatives of $\alpha$ that were needed in the proof of (a) were $\partial_{(x,\theta)}\partial_{(t,r)}\alpha$. 
\end{remark}


%

\begin{thebibliography}{13}


\bibitem{AB}
E. V. Arbuzov,  A. L. Bukhgeim, S. G. Kazantsev.
\newblock
Two-dimensional tomography problems and the theory of $A$-analytic functions.
\newblock {\em Siberian Adv. Math.}, 8(4):1-20, 1998.

\bibitem{Bal-Tamasan}
G.~Bal and A.~Tamasan.
\newblock Inverse source problems in transport equations.
\newblock {\em SIAM J. Math. Anal.}, 39(1):57--76, 2007.

\bibitem{CBG} J. Chang, R. L. Barbour, And  H. Graber. 
\newblock
Imaging of fluorescence in highly
scattering media.
\newblock
{em IEEE Trans. Biomed. Eng.}, 44(9):810-822 1997.

\bibitem{Choulli-Stefanov}
M. Choulli and P. Stefanov.
\newblock Inverse scattering and inverse boundary value problems for the linear
{B}oltzmann equation.
\newblock {\em Comm. Partial Differential Equations}, 21(5-6):763--785, 1996.

\bibitem{COSC}
P. R. Contag, I. N. Olomu, D. K. Stevenson, and C. H. Contag.
\newblock
Bioluminescent
indicators in living mammals.
\newblock
{\em Nature Medicine}, 4(2):245-247, 1998.

\bibitem{FSU}
B.~Frigyik, P.~Stefanov, and G.~Uhlmann.
\newblock The {X}-ray transform for a generic family of curves and weights.
\newblock {\em J. Geom. Anal.}, 18(1):81--97, 2008.

\bibitem{JRB}
H. B. Jang, S. Ramesh, and M. Barlett.
\newblock
Combined optical and fluorescence
imaging for breast cancer detection and diagnosis.
\newblock
{\em Crit. Rev. Biom. Eng.}, 28(3-4):371-375, 
2000.

\bibitem{L}
E. W. Larsen.
\newblock
The inverse source problem in radiative transfer.
\newblock
{\em J. Quant. Spect. Radiat. Transfer},
15:1-5, 1975.

\bibitem{MikhlinP}
S.~G. Mikhlin and S.~Pr{\"o}ssdorf.
\newblock {\em Singular integral operators}.
\newblock Springer-Verlag, Berlin, 1986.
\newblock Translated from the German by Albrecht B\"ottcher and Reinhard
  Lehmann.

\bibitem{Novikov}
R.~G. Novikov.
\newblock An inversion formula for the attenuated {X}-ray transformation.
\newblock {\em Ark. Mat.}, 40(1):145--167, 2002.

\bibitem{P}
A. N. Panchenko.
\newblock
Inverse source problem of radiative transfer: a special case of
the attenuated Radon transform.
\newblock
{\em Inverse Problems}, 9:321-337, 1993.

\bibitem{Reed-Simon1}
M.~Reed and B.~Simon.
\newblock {\em Methods of modern mathematical physics. {I}. {F}unctional
  analysis}.
\newblock Academic Press, New York, 1972.

\bibitem{Si}
C. E. Siewert.
\newblock
An inverse source problem in radiative transfer.
\newblock
{\em J. Quant. Spect.
Radiat. Transfer}, 50:603-609, 1993.

\bibitem{St}
E. Stein.
\newblock {\em Singular integrals and differentiability properties 
of functions}.
\newblock
Princeton University Press, Princeton, New Jersey, 1970.

\bibitem{SU-JAMS}
P.~Stefanov and G.~Uhlmann.
\newblock Boundary rigidity and stability for generic simple metrics.
\newblock {\em J. Amer. Math. Soc.}, 18(4):975--1003, 2005.

\bibitem{Taylor-book1}
M.~E. Taylor.
\newblock {\em Partial differential equations. {I}}, volume 115 of {\em Applied
  Mathematical Sciences}.
\newblock Springer-Verlag, New York, 1996.

\bibitem{YSM}
H. C. Yi, R. Sanchez, and N. J. McCormick.
\newblock Bioluminescence estimation from ocean in situ irradiances.
\newblock
{\em Applied Optics}, 31:822-830, 1992.
\end{thebibliography}

\end{document}